\documentclass[10pt]{article}

\def\today{26~mars~2001}

\begin{document}

\title
{Sur les transformations de Schlesinger de la sixi\`eme \'equation de 
Painlev\'e\footnote{
Note pr\'esent\'ee par Yvonne Choquet-Bruhat,
 remise le 16~octobre~2000, 
           accept\'ee le 12~f\'evrier~2001. \hfill S2001/008.
}
}

\author{Robert Conte
\\
\\ Service de physique de l'\'etat condens\'e, CEA--Saclay
\\ F--91191 Gif-sur-Yvette Cedex, France
\\ Courriel:  Conte@drecam.saclay.cea.fr
}

\maketitle

\noindent \textit{PACS 1995~:}
02.30.Gp, 
02.30.Hq. 
\medskip

\textbf{R\'esum\'e}.
Apr\`es la r\'ecente d\'ecouverte de deux nouvelles transformations de
Schlesinger (TS) pour la sixi\`eme \'equation de Painlev\'e,
nous donnons les relations d'interd\'ependance de toutes les TS connues.
Nous isolons ainsi l'unique d'entre elles qui \`a la fois conserve 
la variable ind\'ependante et n'est pas un produit d'autres TS.
\medskip

\textbf{On the Schlesinger transformations of the sixth Painlev\'e equation}

\textbf{Abstract}.
Following the recent discovery of two new Schlesinger transformations (ST)
for the sixth Painlev\'e equation,
we give the interrelations between all the known STs.
We thus isolate the unique one which at the same time
conserves the independent variable and is not a product of other STs.

\newtheorem{definition}{Definition}
\def\rightnote{\hfill \today. S2001/008}
\def\ccomma{\raise 2pt\hbox{,}} 
\def\AnnENS{Ann.~\'Ec.~Norm.~}
\def\CRAS{C.~R.~Acad.~Sc.~Paris}
\def\LMP{Lett.~Math.~Phys.~}
\def\SAM{Stud.~Appl.~Math.~}

\def\PVI    {{\rm P6}}

\def\Alpha{A}
\def\Beta {B}
\def\abcd{\alpha,\beta,\gamma,\delta}
\def\ABCD{\Alpha,\Beta,\Gamma,\Delta}

\def\veca{{\overrightarrow{\alpha}}}
\def\vecA{{\overrightarrow{\Alpha}}}
\def\vect{{\overrightarrow{\theta}}}
\def\vecT{{\overrightarrow{\Theta}}}

\def\TS{TS}

\def\Sa  {{\rm S}_a}
\def\Sb  {{\rm S}_b}
\def\Sc  {{\rm S}_c}
\def\Sd  {{\rm S}_d}

\def \Pbadc {{\rm H}_{\rm badc}}  
\def \Pdcba {{\rm H}_{\rm dcba}}  
\def \Pcdab {{\rm H}_{\rm cdab}}  

\def \Pca   {{\rm H}_{\rm cbad}}  
\def \Pcb   {{\rm H}_{\rm acbd}}  
\def \Pdb   {{\rm H}_{\rm adcb}}  
\def \Pdc   {{\rm H}_{\rm abdc}}  

\def \TFY   {{\rm T}_{\rm FY}}    
\def \TOk   {{\rm T}_{\rm Ok}}    
\def \TMS   {{\rm T}_{\rm MS}}    
\def \TNJH  {{\rm T}_{\rm NJH}}   
\def \TPVI  {{\rm T}_{\rm CM}}    

\def\Mugan{Mu\u gan}
\section{Introduction}
\label{sectionIntroduction}

Les transformations de Schlesinger (TS) sont d\'efinies \cite{SchlesingerP6}
comme des transformations discr\`etes qui pr\'eservent (\`a des signes pr\`es)
la monodromie d'un syst\`eme diff\'erentiel lin\'eaire appel\'e
syst\`eme de Schlesinger.

Si l'on choisit pour syst\`eme de Schlesinger l'\'equation du second ordre
poss\'edant quatre points fuchsiens (``singuliers-r\'eguliers''),
de birapport $x$,
et une singularit\'e apparente d'affixe $u$,
ses conditions d'isomonodromie \cite{FuchsP6,GarnierThese,JimboMiwaII}
conduisent \`a la contrainte que $u(x)$ doit v\'erifier
l'\'equation ma\^{\i}tresse de Painlev\'e $\PVI$
\label{P6P0def}
\begin{eqnarray*}
\PVI\ : \
u''
&=&
{1 \over 2} \left[{1 \over u} + {1 \over u-1} + {1 \over u-x} \right] u'^2
- \left[{1 \over x} + {1 \over x-1} + {1 \over u-x} \right] u'
\\
& &
+ {u (u-1) (u-x) \over x^2 (x-1)^2} 
  \left[\alpha + \beta {x \over u^2} + \gamma {x-1 \over (u-1)^2} 
        + \delta {x (x-1) \over (u-x)^2} \right]\cdot
\end{eqnarray*}
Aux \TS\ de $\PVI$ d\'ej\`a connues \cite{FY,Okamoto1987I,MS1995b}
viennent de s'en ajouter deux nouvelles \cite{NJH,CM2001OneSemi},
et le but de cette Note est de d\'eterminer la plus \'el\'ementaire de toutes
ces \TS.

Une \TS\ de $\PVI$, qui transforme donc une \'equation
\begin{eqnarray} 
& &
E(u) \equiv \PVI(u,x,\veca)=0,\
\veca=(\abcd),\
\label{eqPnu}
\end{eqnarray}
en une autre \'equation de m\^eme forme
\begin{eqnarray} 
& &
E(U) \equiv \PVI(U,X,\vecA)=0,\
\vecA=(\ABCD),
\label{eqPnU}
\end{eqnarray}
est repr\'esentable d'au moins trois mani\`eres~:
\begin{enumerate}
\item
dans l'espace $(u,x,U,X)$,
par une \textit{transformation birationnelle}
($u$ rationnel en $(U,U')$, $U$ rationnel en $(u,u')$),
accompagn\'ee d'une \'eventuelle homographie entre $x$ et $X$,

\item
dans l'espace des param\`etres $(\veca,\vecA)$,
par quatre relations alg\'ebriques,

\item
dans l'espace des exposants de monodromie $(\vect,\vecT)$,
d\'efinis par les relations
\begin{eqnarray} 
& &
\theta_{\infty}^2= 2 \alpha,\
\theta_0^2       =-2 \beta,\
\theta_1^2       = 2 \gamma,\
\theta_x^2       =1 - 2 \delta,
\end{eqnarray}
par une transformation affine
\begin{eqnarray}
& &
\vect = M_1 \vecT + M_0,\
\vect=\pmatrix{\theta_{\infty} \cr \theta_0 \cr \theta_1 \cr \theta_x},\
\vecT=\pmatrix{\Theta_{\infty} \cr \Theta_0 \cr \Theta_1 \cr \Theta_x},\
\label{eqSTM1M0}
\end{eqnarray}
o\`u $M_1$ et $M_0$ sont des matrices de nombres rationnels.

\end{enumerate}

\section{Les transformations de Schlesinger connues}

Nous les \'enum\'erons dans l'ordre chronologique\footnote
{
Apr\`es soumission,
nous avons d\'ecouvert l'ant\'eriorit\'e de Garnier \cite{Garnier1943a}
sur les R\'ef.~\cite{FY,Okamoto1987I,MS1995b}.
Il a \'etabli une formule \'equivalente \`a  (\ref{eqTMSu}).
}.

En appliquant deux fois une transformation alg\'ebrique
entre $\PVI$ et l'\'equation d'ordre deux et de degr\'e deux
satisfaite par $(q-1) (q-x) p$,
o\`u $q=u$ et $p$ sont la position et l'impulsion de l'hamiltonien de $\PVI$
donn\'e par Malmquist \cite{MalmquistP6},
Fokas et Yortsos \cite{FY} ont obtenu une transformation entre deux 
\'equations $\PVI$ repr\'esent\'ee par
\begin{eqnarray}
& &
\TFY : \vect = \pmatrix{0 & 1 & 0 & 0 \cr
                        1 & 0 & 0 & 0 \cr
                        0 & 0 & 0 & 1 \cr
                        0 & 0 & 1 & 0 \cr} \vecT
              + \pmatrix{ 1 \cr 1 \cr 0 \cr 0 \cr},\
x=\frac{1}{X}\cdot
\label{eqFY}
\end{eqnarray}
La transformation entre $u$ et $U$,
produit de deux transformations non-birationnelles,
est bien birationnelle,
comme le montrent les relations (\ref{eqTMSu})
et (\ref{eqFYOkMS1})--(\ref{eqFYOkMS2}) ci-dessous.

Par l'\'etude des transformations canoniques de l'hamiltonien de $\PVI$
\cite{MalmquistP6},
Okamoto \cite{Okamoto1987I} a obtenu la \TS\ suivante,
canonique au sens hamiltonien,
\begin{eqnarray}
& &
\TOk : \vect = \vecT + \pmatrix{ 1 \cr 0 \cr 0 \cr 1 \cr},\
x=X.
\label{eqOkamoto}
\end{eqnarray}
Son \'ecriture birationnelle 
r\'esulte des relations (\ref{eqTMSu})
et (\ref{eqFYOkMS1})--(\ref{eqFYOkMS2}) ci-dessous.

En partant d'une paire de Lax matricielle obtenue par isomonodromie
\cite{JimboMiwaII},
et apr\`es r\'esolution d'un probl\`eme de Riemann-Hilbert,
\Mugan\ et Sakka \cite{MS1995b} ont obtenu une \TS\ 
que nous pouvons \'ecrire sous la forme factoris\'ee
\begin{eqnarray}
& &
{\hskip -1.0 truecm}
\frac{(u-x)(U-x)}{(x-1) u U} = \frac{R_n^{+} R_n^{-}}{R_d^{+} R_d^{-}}\ccomma\
x=X,\
\label{eqTMSu}
\\
& &
{\hskip -1.0 truecm}
R_n^{\pm} =x (1-x) U' 
          + \theta_{\infty} (U-1) (U-x)
          - \Theta_{1} (U-x)
  + (1 \pm \theta_{x}) x (U-1),
\\
& &
{\hskip -1.0 truecm}
R_d^{\pm} =x (1-x) U' 
          + \theta_{\infty} U (U-1)
          + \Theta_{1}      (x-1) U
        \pm \theta_{0}     x (U-1),
\\
& &
{\hskip -1.0 truecm}
\TMS : \vect = \vecT + \pmatrix{1 \cr 0 \cr 1 \cr 0 \cr}\cdot
\label{eqTMS}
\end{eqnarray}

En exploitant les propri\'et\'es d'une \'equation aux d\'eriv\'ees partielles
de type ``schwarzien'',
Nijhoff \textit{et al.} \cite{NJH}
ont trouv\'e une \TS\ represent\'ee par
\begin{eqnarray}
& &
(u-x) U -1 
=
\frac
{2 (\Theta_{\infty} - \theta_x) (u-1) (u-x)}
{x (1-x) u' 
+ \theta_{\infty} (u-1)(u-x)
- \theta_1 (u-x)
+ (\theta_x +1) x (u-1)}\ccomma
\\
& &
\TNJH: \vect = \frac{1}{2} \pmatrix{-1 &  1 & -1 & -1 \cr
                                    -1 & -1 &  1 & -1 \cr
                                    -1 & -1 & -1 &  1 \cr
                                     1 & -1 & -1 & -1 \cr} \vecT
             + \frac{1}{2} \pmatrix{ 1 \cr 1 \cr 1 \cr 1 \cr},\
x-1=\frac{1}{X-1}\cdot
\label{eqNJH}
\end{eqnarray}

En \'etendant une m\'ethode de ``troncature'',
nous avons obtenu la \TS\ \cite{CM2001OneSemi}
\begin{eqnarray}
\frac{2(\theta_\infty-\Theta_\infty)}{u-U}
& = &
\frac{x(x-1)U'}{U(U-1)(U-x)}
+\frac{\Theta_0}{U}
+\frac{\Theta_1}{U-1}
+\frac{\Theta_x-1}{U-x}
\label{eqTP6uvect}
\\
& = &
\frac{x(x-1)u'}{u(u-1)(u-x)}
+\frac{\theta_0}{u}
+\frac{\theta_1}{u-1}
+\frac{\theta_x-1}{u-x}
\label{eqTP6Uvect}
\\
\TPVI & : &
        \vect = \frac{1}{2} \pmatrix{-1 &  1 & -1 & -1 \cr
                                      1 & -1 & -1 & -1 \cr
                                     -1 & -1 & -1 &  1 \cr
                                     -1 & -1 &  1 & -1 \cr} \vecT
              + \frac{1}{2} \pmatrix{ 1 \cr 1 \cr 1 \cr 1 \cr},\
x=X.
\label{eqT6}
\end{eqnarray}

La  transformation de Kitaev \cite{Kitaev1991}
ne figure pas dans la liste ci-dessus.
C'est une application \`a $\PVI$,
suppos\'ee \'ecrite avec une somme de quatre fonctions elliptiques
\cite{PaiCRAS1906},
de la transformation de Landen entre fonctions elliptiques
\cite{Manin1998}.
Cependant,
elle n'existe qu'au prix de deux contraintes entre les quatre param\`etres de
$\PVI$
\begin{eqnarray}
& &
\theta_0-1=\theta_x,\
\theta_1=\theta_{\infty},
\end{eqnarray}
ce qui lui interdit de d\'efinir une \TS.

\section{La transformation de Schlesinger \'el\'ementaire de $\PVI$}
\label{sectionElementaryPVI}

Pour comparer entre elles ces diverses TS,
nous utilisons la repr\'esentation affine (\ref{eqSTM1M0}),
de beaucoup la plus simple. Elle est donn\'ee par les cinq relations
(\ref{eqFY}), 
(\ref{eqOkamoto}),
(\ref{eqTMS}), 
(\ref{eqNJH}) et 
(\ref{eqT6}).

Deux classes de transformations banales laissent $\PVI$ invariante.
Il s'agit d'une part des quatre changements ind\'ependants de signes de
$\theta_{\infty},\theta_{0},\theta_{1},\theta_{x}$,
que nous notons $\Sa,\Sb,\Sc,\Sd$,
par exemple pour $\theta_{\infty}$,
\begin{eqnarray}
& &
\Sa : \vect = 
 \pmatrix{-1 & 0 & 0 & 0 \cr
           0 & 1 & 0 & 0 \cr
           0 & 0 & 1 & 0 \cr
           0 & 0 & 0 & 1 \cr} \vecT,\
\theta_{\infty}=-\Theta_{\infty},\ 
x=X,\
u=U,\ 
\end{eqnarray}
et d'autre part des $24$ \'el\'ements du groupe des permutations
des quatre points singuliers de $\PVI$,
represent\'es par des homographies simultan\'ees sur $u$ et sur $x$.
Le sous-groupe qui conserve $x$ est constitu\'e de quatre \'el\'ements,
l'identit\'e et les trois homographies
\begin{eqnarray}
& &
\Pbadc : \vect = 
 \pmatrix{ 0 & 1 & 0 & 0 \cr
           1 & 0 & 0 & 0 \cr
           0 & 0 & 0 & 1 \cr
           0 & 0 & 1 & 0 \cr} \vecT,\
x=X,\
u=\frac{x}{U}\ccomma
\\
& &
\Pdcba : \vect = 
 \pmatrix{ 0 & 0 & 0 & 1 \cr
           0 & 0 & 1 & 0 \cr
           0 & 1 & 0 & 0 \cr
           1 & 0 & 0 & 0 \cr} \vecT,\
x=X,\
u-x=\frac{x(x-1)}{U-x}\ccomma
\\
& &
\Pcdab : \vect = 
 \pmatrix{ 0 & 0 & 1 & 0 \cr
           0 & 0 & 0 & 1 \cr
           1 & 0 & 0 & 0 \cr
           0 & 1 & 0 & 0 \cr} \vecT,\
x=X,\
u-1=\frac{1-x}{U-1}\ccomma
\\
& &
\Pbadc^2=\Pdcba^2=\Pcdab^2=1,\
\Pbadc   \Pdcba   \Pcdab = 1.
\end{eqnarray}

Les diverses {\TS} ob\'eissent aux relations
\begin{eqnarray}
& &
\TFY = \Pbadc \Pcb \TMS \Pcb,\
\TMS = \Pcb \Pbadc \TFY \Pcb,\
\label{eqFYOkMS1}
\\
& &
\TOk = \Pdc \TMS \Pdc,\
\TMS = \Pdc \TOk \Pdc,\
\label{eqFYOkMS2}
\\
& &
\TNJH=\Pdb \TPVI \Pbadc,\
\TPVI=\Pdb \TNJH \Pbadc, 
\label{eqNJHPVI}
\\
& &
\TOk=\Pdcba (\TPVI \Pbadc \Sa \Sd)^2 
    =\Pca \Sa \Sb   (\Sa     \TPVI \Pbadc)^3 \Sc \Sd \Pca,  
\label{eqOkPVI}
\\
& &
\TPVI^2=1,\
(\Sa \TPVI)^3=1,\
(\Sa \Sb \TPVI \Pbadc)^4=1,\
(\Sa \TPVI \Pbadc)^6=1,\     
\label{eqPVIpower}
\end{eqnarray}
 o\`u les diverses transpositions 
agissent sur $(u,x)$ par
\def \Pca   {{\rm H}_{\rm cbad}}  
\def \Pcb   {{\rm H}_{\rm acbd}}  
\def \Pdb   {{\rm H}_{\rm adcb}}  
\def \Pdc   {{\rm H}_{\rm abdc}}  
\begin{eqnarray}
& &
 \Pdb\ : \ x-1=\frac{1}{X-1}\ccomma\ u-x = (1-x) U,  
\\
& &
 \Pca\ : \ x-1=\frac{1}{X-1}\ccomma\ u-1 = \frac{1}{U-1}\ccomma  
\\
& &
 \Pdc\ : \ x=1/X,\ u=x U,                            
\\
& &
 \Pcb\ : \ x=1-X,\ u=1-U.                            
\end{eqnarray}
Ces relations prouvent
\begin{enumerate}
\item
 l'\'equivalence entre $\TFY$, $\TOk$ et $\TMS$ 
 (relations (\ref{eqFYOkMS1})--(\ref{eqFYOkMS2})),
\item
 l'\'equivalence entre $\TNJH$ et $\TPVI$ (relation (\ref{eqNJHPVI})),
\item
 la non-\'el\'ementarit\'e de $\TFY$, $\TOk$, $\TMS$ 
(relation (\ref{eqOkPVI}))
 qui, \`a des signes et des homographies pr\`es,
sont le carr\'e ou le cube de $\TPVI$.
\end{enumerate}

Enfin,
il importe d'essayer de conserver aussi la variable ind\'ependante $x$,
c'est par exemple indispensable pour \'etablir une relation de 
contigu\"{\i}t\'e \cite{CM2001OneSemi}.
L'unique \TS\ qui soit \'el\'ementaire et qui conserve $x$ est alors $\TPVI$.
\`A des signes et des homographies pr\`es,
c'est elle-m\^eme une racine de l'unit\'e
(relation (\ref{eqPVIpower})),
d'ordre deux, trois, quatre ou six.

\medskip
\textbf{Remerciements}.
L'auteur remercie M.~Musette pour de fructueuses discussions.


\vfill \eject
\end{document}